# A Fast Relax-and-Round Approach to Unit Commitment for Data Center Own Generation

Shaked Regev, Eve Tsybina, Slaven Peles
Oak Ridge National Laboratory, Oak Ridge, TN, USA
Email: {regevs, tsybinae, peless} @ ornl.gov

***Abstract* — The rapid growth of data centers increasingly requires data center operators to "bring own generation" to complement the available utility power plants to supply all or part of data center load. This practice sharply increases the number of generators on the bulk power system and shifts operational focus toward fuel costs rather than traditional startup and runtime constraints. Conventional mixed-integer unit commitment formulations are not well suited for systems with thousands of flexible, fast-cycling units. We propose a unit commitment formulation that relaxes binary commitment decisions by allowing generators to be fractionally on, enabling the use of algorithms for continuous solvers. We then use a rounding approach to get a feasible unit commitment. For a 276-unit system, solution time decreases from 10 hours to less than a second, with no accuracy degradation. Our approach scales with no issues to tens of thousands of generators, which allows solving problems on the scale of the major North America interconnections. The bulk of computation is parallel and GPU compatible, enabling further acceleration in future work.**

## I. Introduction

DATA center expansion is exerting increasing pressure on power systems. This pressure arises not only from new load but also from increased system complexity as data-center developers procure and interconnect generation resources. According to 2024 S&P Global Market Intelligence report [1], in the past five years datacenters have added up to 27.8 GW of utility power load in the US and almost 1 GW of utility power load in Canada. An additional 11.6 GW are under construction in the US as of Q4 2024 and 43.8 GW are planned for deployment during 2025-2030. As of today, most of that load is being developed in a few states, such as Virginia, Texas, and increasingly, Arizona and Georgia, but the geography of data centers is likely to increase as more regions announce data center development plans [2]. The projected data centers are larger, and more complex in nature than those in operation, with more advanced hardware architectures and higher energy intensity per rack [3].

This growth in data center development has resulted in electricity supply being unable to keep up with the growing demand. As has been witnessed by a number of recent regulatory and industry meetings [4]-[5], data centers are increasingly expected to "bring own generation". This term refers to investor-owned generation assets intended to maintain resource adequacy of transmission operators. Further, several large investors have begun acquiring utility generation fleets to delay planned retirements and support accelerated data-center deployment schedules [6]. The data center investment in power generation is likely to result in a very characteristic generation profile: the new units are likely to be small- and medium-sized, fast ramping, and to have higher fuel costs due to the smaller scale and the use of more expensive fuels such as natural gas. Further, the owners of such generating assets can be expected to prefer partial load of their units in order to maintain spinning reserves for higher reliability of grid operation.

Meanwhile, the present best practices for unit commitment are not fit for addressing problems with a large number of flexible resources, as they are based on solving mixed integer programs (MIPs) [7]. The MIP problem is NP-hard because an exact solution requires checking all value combinations of the integer variables. To address this challenge, current state-of-the-art MIP methods, such as those used in Gurobi [8], CPLEX [9] and other commercial tools, use heuristics, such as branch-and-bound [10]-[11]. However, even with those heuristics, MIP methods still have poor scaling properties and become inefficient for large problems [11]. We attempt to address the need for a faster and more efficient relax-and-round unit commitment (RRUC) process by relaxing the integer problem and creating a new algorithm that allows to commit thousands of generating units in a matter of minutes.

The study is organized as follows. Section II describes the methodology and solver inputs. Section III presents the numerical results. Section IV discusses implications and outlines directions for future research.

Notice: This manuscript has been authored by UT-Battelle, LLC, under contract DE-AC05-00OR22725 with the US Department of Energy (DOE). The US government retains and the publisher, by accepting the article for publication, acknowledges that the US government retains a nonexclusive, paid-up, irrevocable, worldwide license to publish or reproduce the published form of this manuscript, or allow others to do so, for US government purposes. DOE will provide public access to these results of federally sponsored research in accordance with the DOE Public Access Plan (https://www.energy.gov/doe-public-access-plan). Research sponsored by the Laboratory Directed Research and Development Program of Oak Ridge National Laboratory, managed by UT-Battelle, LLC, for the US Department of Energy.

## II. Mathematical Formulation and Methodology

### A. *Mathematical formulation*

Recall, generating units brought by data center operators are likely to be different from the existing utility generating units, based on their primary operational role. They will be predominantly small- to medium-sized, with fast ramping capability, most likely small fossil fuel units. Importantly, a portion of this capacity is likely to be maintained as spinning reserve rather than dispatched at full output, to provide data center reliability. The associated strategic behavior in power markets is usually characterized by very high prices in the right half of the supply curves ("tail" effect).

Finally, data center loading conditions are relatively new to power systems research, and may not be well known. Preliminary evidence from the ORNL Frontier supercomputer suggests that about 95% of total power consumed by a supercomputing facility is active power for the GPU stack and overhead, with about 5% for reactive power cooling system needs. This allows us to assume that modeling only active power constitutes a reasonable first step in data center dominated power systems.

Combined, these characteristics indicate that the "own generation" units will have highly convex electricity market bids and be capable of frequent cycling and short startup times. This allows for the preliminary formulation of an objective function that emphasizes cost structure compared to startup and runtime constraints. This unit commitment formulation is specified in Eq. 1.

$$\min_{P_i} \sum_{i=1}^{n} u_i(a_i P_i^2 + b_i P_i + c_i) \tag{1}$$

s.t. $P_{min,i} \le P_i \le P_{max,i}$, $u_i \in \{0,1\}$
$\sum_{i=1}^{n} u_i P_i \ge D$,
$\sum_{i=1}^{n} u_i\, P_{max,i} \ge D + 3\sigma_D + \max_l P_{max,l}$

where $u_i$ is a binary variable that indicates whether generator $i$ is on or off, $P_i$ is active power of the generating unit $i$, which is supplied if a unit gets dispatched, D is the total power demand, $\sigma_D$ is a standard deviation of D, and $a_i, b_i \ge 0$ and $c_i$ are constant parameters. The last constraint is what limits the risk of the solution, by ensuring that enough generators are committed in the case of a contingency, or if demand is higher than expected. The active power is bound by the minimum stable output $P_{min,i}$ and maximum capacity $P_{max,i}$, usually set in accordance with technological limitations of the generating equipment. The sum of total committed capacity should be larger than or equal to the total system load. For the purposes of numerical simulations, we set load to be equal to 80% of the preferred committed capacity. The problem is specified to be convex, therefore the quadratic cost coefficient $a_i$ is nonnegative. The cost function is nondecreasing, therefore the linear coefficient $b_i$ is nonnegative.

Since we are primarily interested in the commitment of fast small generating units, we do not incorporate unit commitment components such as intraday constraints, maximum number of startups, or runtime constraints. The main questions addressed by the study relate to computational scalability, effects of high fuel costs, and effects of strategically running turbines below full capacity.

Our approach starts by relaxing discrete variables$u_i \in 0,1$ to continuous ones $y_i \in [0,1]$ and performing continuous nonlinear constrained optimization of the relaxed problem: $y_i(a_i P_i^2 + b_i P_i + c_i)$. This step is similar to the branch-and-bound method. Such analysis will typically produce non-integer values for many, if not all, of $y_i$.

However, instead of branching for $y_i$ with non-integer values, we reorder terms in the descending order in $y_j$ so that $y_j > y_{j+1}$. Then, since in our problem the $y_j$ variables represent which generators are committed and which are not, we find the first $m$ generators such that $\sum_{j=1}^{m} P_{max,j} \ge D + 3\sigma_D + \max_j P_{max,j}$. We then solve economic dispatch (Eq. (2)) with all options of generators that can supply the load $k{:}\ s.t.\ m \le k \le n$, (meaning we set $y_{j\le k} = 1$ and $y_{j>k} = 0$).

$$\min_{P_j} \sum_{j=1}^{k} \left(a_j P_j^2 + b_j P_j + c_j\right) \tag{2}$$

s.t. $P_{min,j} \le P_j \le P_{max,j}$
$\sum_{j=1}^{k} P_j \ge D$,
$\sum_{j=1}^{k} P_{max,j} \ge D + 3\sigma_D + \max_l P_{max,l}$

Testing $\forall k{:}\ s.t.\ m \le k \le n$ requires solving $n - m + 1$ completely independent instances of Eq. (2), this can be done in parallel. Solving Eq. (1) gives us a preference of the generators, which reduces our search space from exponential in the number of generators to linear. This makes brute forcing over the remaining solutions tractable. Other implementations could check only some values of $k$.

The time complexity of this approach is that of $n - m + 2$ continuous constrained nonlinear optimization runs. The first run solves Eq. (1) and the other $n - m + 1$ solve Eq. (2), possible in "embarrassingly parallel" fashion. The complexity of the single run is dominated by the complexity of the sparse linear solver used underneath, which is typically $\approx O(N^{1.5})$, where $N$ is the size of the optimization problem. The time complexity of a branch-and-bound approach without pruning is $\approx O(N^{1.5} b^d)$, where $b$ is the maximum number of children at each branching point, and $d$ is the search depth of the tree [11]. Even with pruning rules, branch-and-bound scales poorly, as our experiments in Section III.A show.

We point out that in this work we use a simple ordering in $y_j$. However, our approach allows for seamless inclusion of low-complexity algorithms that would improve the solution by incorporating more domain information, as demonstrated for example in [12].

### B. Generator selection

To verify the methods proposed above, we simulate a lossless DCOPF system with an arbitrary number of generators, each having its own unit behavior. The startup, shutdown, and operating cost components are modelled based on real PJM generators, available from PJM Data Miner [12]. The original dataset for 12 months of PJM market bids contains 488 unique generating units that have the minimum runtime of less than or equal to 1 hour, and among those 46 units that have a sufficient number of steps in supply curves to approximate a cost function. The 46 units submitted approximately 330,000 bid curves throughout 2024, including 18,677 bid curves that were unique. The properties of the selected bid curves are summarized using Fig. 1. The units tend to have $P_{min}$ of about 20-24 MW, Pmax of up to 50-100 MW, with only a few units as large as 200-300 MW. The no load cost tends to be over 200 USD/hour. The startup cost for many units is comparatively small, up to 100 USD per start.

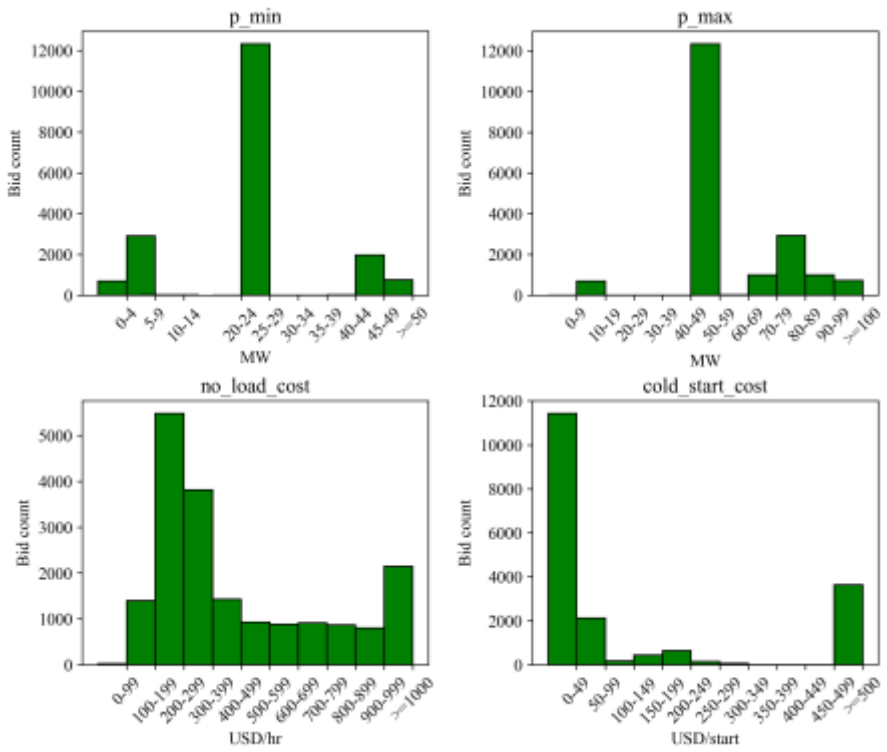

**Fig. 1.** Properties of the selected bid curves for PJM generating units.

The variable cost properties are illustrated in Fig. 2. Each of the 46 unique generators has 3-7 steps in supply curves, and some units have a visible "tail" indicating that a unit owner prefers to not be dispatched for the volume in higher steps of the bid curve. The PJM supply curves are not assumed to exactly represent marginal costs of production due to strategic bidding behaviors. However, we believe that they are representative of realistic bids that transmission operators can expect from the additional generating units.

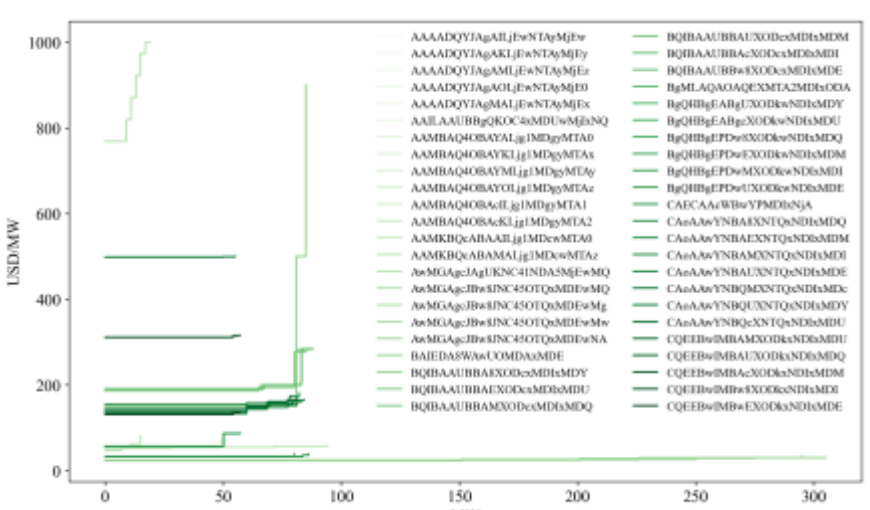

**Fig. 2.** Snapshots of supply curves for the 46 unique generating units.

To convert the information about bids into quadratic total costs per MWh, we add no load costs and, where applicable, startup costs as fixed cost component. To obtain the curve properties for variable costs, we first integrate the marginal cost bid curves in Figure 2 to obtain the total cost. Then, we fit quadratic cost curves to the total cost using least squares. In all cases $R^2 > 0.998$, though in some $a \approx 0$. This indicates that a quadratic curve is a reasonable of operation costs, though sometimes linear is sufficient.

## III. Simulation and Results

### A. Computational scalability

We use the quadratic functions, startup, and no-load costs obtained in the previous step of the study to simulate the commitment of varying number of units. We implement our method using Julia programming model [14] and MadNLP optimization package [15], which is based on our prior work [16]-[17]. As a reference MIP implementation, we choose the open-source Juniper package [18], which runs on the same JuMP platform [19] as MadNLP.

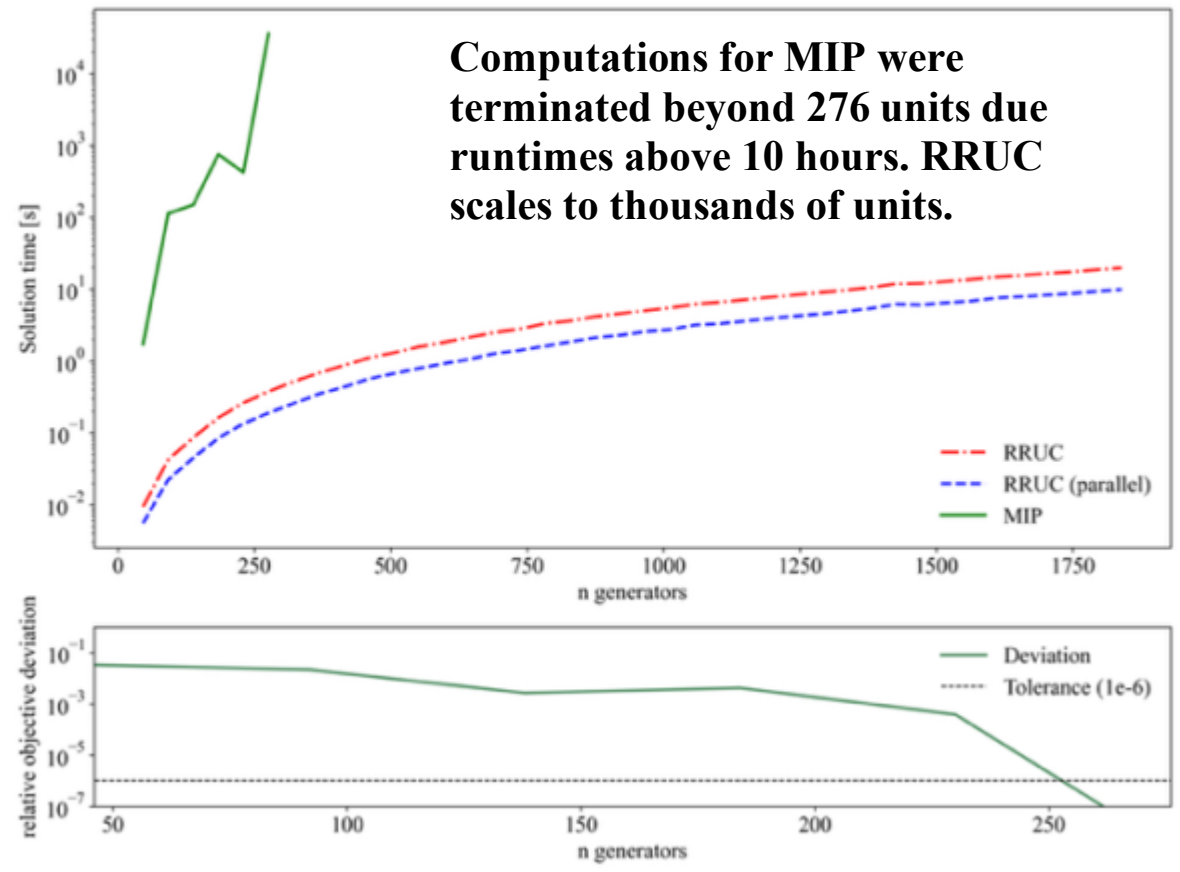

**Fig. 3.** Test results for 46-1840 generating units, (a) solve time on a logarithmic scale, (b) deviation in objective value. For larger problems, the deviation is below the tolerance set for the solvers. All tests used Julia @benchmark, to ensure that both scripts generate consistent solve times.

To assess how the algorithm scales compared MIP, we increased the system size in multiples of 46, by replicating the original sample of 46 unique generating units, with small deviations to their parameters. With a multiplier of 6, i.e. 276 dispatched units, RRUC performed thousands of times faster than the MIP algorithm. We also found its timing to be more consistent. Because MIP problems depend on branch structure, for an irregular tree the algorithm can have varying depths of branches. By contrast, the relaxed problem does not rely on tree searches, therefore it has no risk of having unusually shallow or deep structures. This results in much lower time complexity, as well as more predictable convergence time that scaled monotonously with the system size. Figure 3(a) shows the results of the simulation for up to 1840 units. The MIP solution time explodes, while RRUC scales very well.

## B. *Effects on total cost of the objective function*

In addition to the difference in computational efficiency, we find a difference in value of the objective function for the system. Figure 3(b) shows the deviation of the RRUC objective from the MIP objective, calculated as (RRUC-MIP)/MIP. This number is positive and stays within 6%. For larger problems, it is within the solver tolerance. This is because the MIP solver must also perform approximations, which it does not need to do for smaller problems. To understand possible cost differences, we do a unit-by-unit review of the sets of committed equipment. Each quadratically approximated generating unit has a mapping to the piecewise linear generating unit from PJM data. This allows comparing solutions based on real supply functions, $P_{min} \doteq \sum_{i=1}^{n} u_i P_{min,i}$, $P_{max,i} \doteq \sum_{i=1}^{n} u_i P_{max,i}$. Both algorithms do not guarantee a global minimum, therefore tracing the exact composition of committed equipment does not guarantee the best explanation of the different objective functions. However, it offers some preliminary evidence.

For the original problem, the algorithms committed 41 units each and only one different unit. For MIP, this was the smallest peak unit, which supplies between 7 and 20 MW at the starting price of 768 USD/MWh, and a cold start cost of USD 227 (upper left in Figure 2). RRUC committed a cheaper but larger unit, from the middle of the supply stack. That unit supplies between 50 and 83 MW at a starting price of 143 USD/MWh, and a cold start cost of USD 5066. MIP's $P_{min}$ was 1558 MW, $P_{max}$ was 3611 MW. RRUC's $P_{min}$ was 1601 MW and $P_{max}$ was 3675 MW. For the scaled problems with 92-230 units, RRUC's $P_{min}$ and $P_{max}$ vary slightly (usually larger) from that of MIP. For 276 units, the same generators are picked and the solutions are identical to within our tolerance. Figure 4 shows the relative difference, computed as (RRUC-MIP)/MIP for $P_{min}$ and $P_{max}$. Preliminary comparison suggests that with few generators, MIP produces smaller system sizes compared to RRUC, and RRUC potentially systematically picks larger units with lower fuel costs and higher startup costs. With more generators, MIP must approximate and becomes almost identical to RRUC.

The tradeoff between low startup costs and high operating costs may present challenges for multi-hour unit commitment problems, as was discussed above. This is left for further research. Further, the computational challenges associated with MIP, combined with the absence of the global optimum in both algorithms, require extensive dedicated effort for the system sizing simulations. This may prove computationally prohibitively expensive. It is possible to use additional data for RRUC to investigate the algorithm's behavior, but this is beyond the scope of the paper.

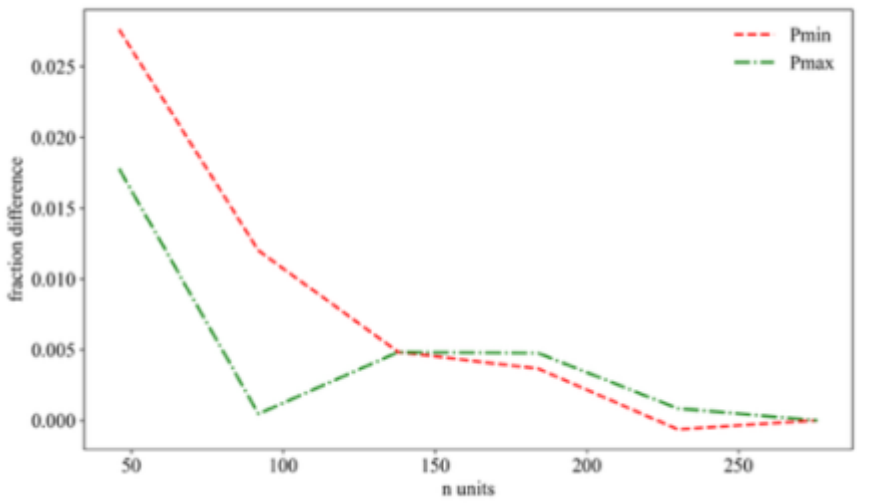

**Fig. 4.** Difference in $P_{min}$ and $P_{max}$ between MIP and RRUC.

## C. *Effects of fixed and variable costs in RRUC*

The fast runtime of RRUC allows us to run very large unit commitment problems and identify fine statistical details about properties of the committed sets of generating units. We are mainly interested in how fixed and variable costs affect unit commitment and system $P_{min}$.

We use the 18,677 bid unique curves that were provided to PJM in 2024 by the 46 studied units. Since those supply curves relate to the same equipment but reflect slight changes in fuel costs and bidding strategies, we can use them as independent supply curves for technologically similar generating units. $P_{min}$ of the proposed large test system is 523.8 GW and $P_{max}$ of the system is 1310.3 GW. We simulate load increments of 100 GW to generate systems that would require committing an increasing number of units, as more favorable units are depleted. Lower steps, e.g. load of 100 GW or 300 GW, mean that the algorithm has a lot of freedom to pick optimal units. The difference between committed units and not committed units would be higher, accurately showing the features that the algorithm prefers. As the load increases to about 800-1000 GW, the algorithm will have less choice, causing the feature values to converge.

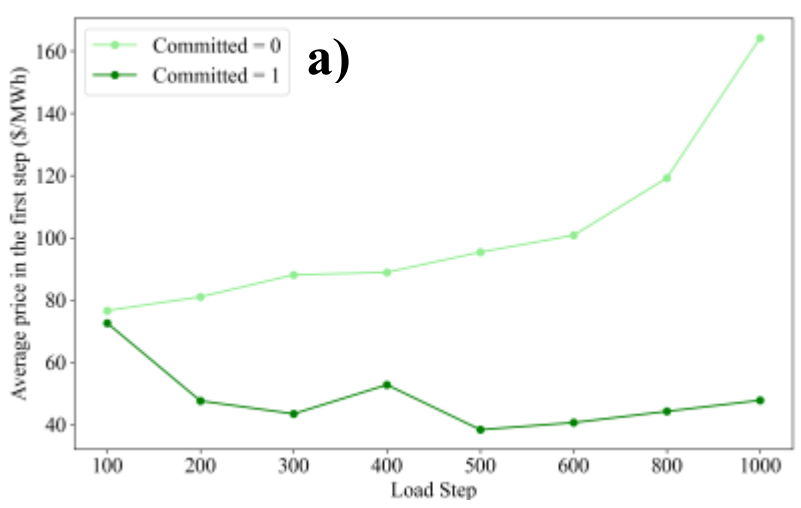

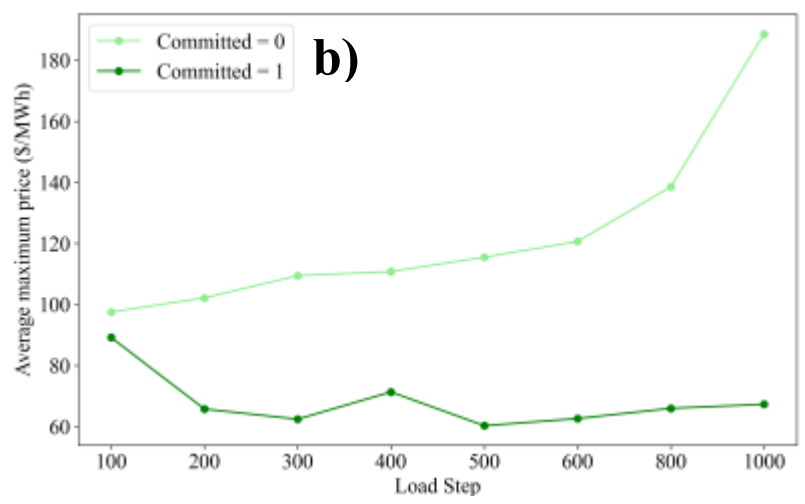

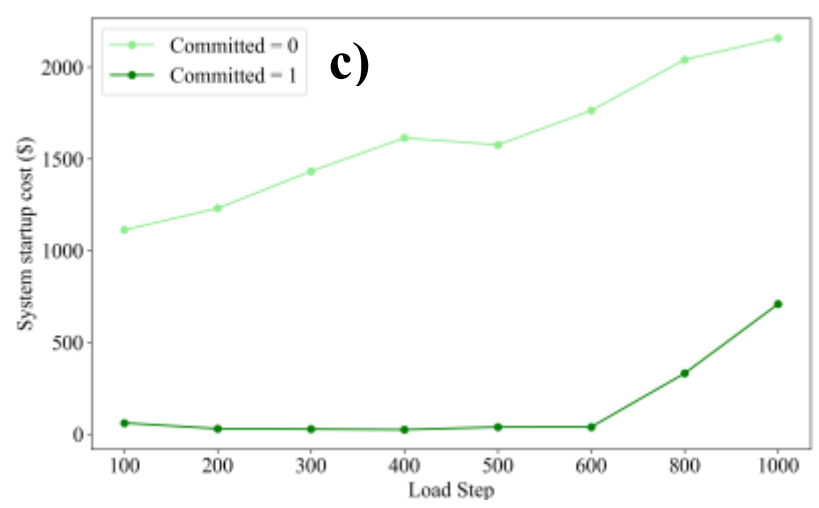

**Fig. 5.** Test results for different levels of freedom in selecting the generating units for RRUC: (a) average price in the first step of the supply curve, (b) average maximum price in the supply curve, (c) total startup costs of committed versus remaining units.

The results of the numerical tests are shown in Figure 5. The algorithm tends to select the generators with lower startup costs, even at the expense of the comparatively high prices in supply curve. But as generators with lowest startup costs get depleted by load step of 600 GW, it increasingly uses units with lower prices in the supply curve and becomes less sensitive to fixed costs. In the input data, the lower bid prices were associated with higher $P_{min}$, which causes the algorithm to commit additional capacity to benefit from the lower prices in the supply curve.

## IV. Conclusions and future work

This study proposes a new methodology for unit commitment, considering the technological and strategic characteristics of the new generation which are likely to come to the grid due to data center expansion. Specifically, we prioritize scalability and the ability of the methods to accommodate small and medium sized generators with high fuel cost components and fewer runtime constraints.

We performed a series of tests to identify the computational and cost performance of our algorithm, and understand the possible reasons behind such performance. Specifically, we find that the proposed algorithm has very high scaling ability. It takes only 0.19 seconds to run 276 units for RRUC but about 10 hours on MIP. Further scaling of the system for algorithm comparison purposes turned to be prohibitively computationally expensive due to excessive runtime of MIP. RRUC also tends to select more expensive composition of generating units, potentially due to selecting larger units with lower fuel costs and higher startup costs. The relative difference between RRUC and MIP were within 6% for the smaller tested system sizes. The relative difference decreases with the problem size and is smaller than our tolerance for the largest system MIP could solve.

An advantage of our approach is that it can be easily augmented by adding additional criteria for ordering generators to commit. For example, in [12] authors demonstrated how low complexity algorithms could be used to enhance unit commitment analysis.

Future work can be divided among improved implementation, adding further constraints such as runtime or must run units, and transmission constraints. Our implementation can further be improved by parallelizing the more substantial part of the computation and moving the solver to the GPU.

## References

[1] 2024 S&P Global Market Intelligence Datacenters and energy report (subscription based).

[2] New Data Center Developments: September 2025, article in Data Center Knowledge Resource Library, https://www.datacenterknowledge.com/data-center-construction/new-data-center-developments-september-2025, accessed 11/14/2025.

[3] Key Trends and Technologies Impacting Data Centers in 2024 and Beyond, article in Data Center Knowledge Resource Library, https://www.datacenterknowledge.com/next-gen-data-centers/key-trends-and-technologies-impacting-data-centers-in-2024-and-beyond, accessed 11/14/2025.

[4] FERC Technical Conference Regarding the Challenge of Resource Adequacy in RTO and ISO Regions, June 4-5, 2025

[5] Shapiro among governors whose states rely on PJM who want data centers to guarantee their own power, article in The Daily Review, https://www.thedailyreview.com/ap/state/shapiro-among-governors-whose-states-rely-on-pjm-who-want-data-centers-to-guarantee-their/article_09f35af3-559c-46d3-af11-5fd531fb7b03.html, accessed 11/14/2025.

[6] Tech Firms Are Building Their Own Power Plants, article in National Association of Manufacturers Input Stories, https://nam.org/tech-firms-are-building-their-own-power-plants-34960/, accessed 11/14/2025.

[7] Achterberg, T., Wunderling, R. (2013). “Mixed Integer Programming: Analyzing 12 Years of Progress”. In: Jünger, M., Reinelt, G. (eds) Facets of Combinatorial Optimization. Springer, Berlin, Heidelberg. https://doi.org/10.1007/978-3-642-38189-8_18.

[8] Gurobi Optimization, LLC. “Gurobi Optimizer Reference Manual”. 2024. https://www.gurobi.com

[9] IBM (2017) IBM ILOG CPLEX 12.7 User’s Manual (IBM ILOG CPLEX Division, Incline Village, NV).

[10] E. L. Lawler, D. E. Wood, (1966) Branch-and-Bound Methods: A Survey. Operations Research 14(4):699-719.

[11] Morrison, D. R., Jacobson, S. H., Sauppe, J. J., and Sewell, E. C. "Branch-and-bound algorithms: A survey of recent advances in searching, branching, and pruning." Discrete Optimization 19 (2016): 79-102.

[12] Holt, Timothy, Abhyankar, Shrirang, Kuruganti, Teja, Schenk, Olaf and Peles, Slaven. “Data-Driven Unit Commitment Refinement-a Scalable Approach for Complex Modern Power Grids”. Proceedings of the 57th Hawaii International Conference on System Sciences, Electric Energy Systems, p. 3082-3092, 2024.

[13] PJM Data Miner, Energy Market Generation Offers.

[14] Bezanson, Jeff, Alan Edelman, Stefan Karpinski, and Viral B. Shah. "Julia: A fresh approach to numerical computing." SIAM review 59, no. 1 (2017): 65-98.

[15] Pacaud, François, and Sungho Shin. "GPU-accelerated dynamic nonlinear optimization with ExaModels and MadNLP." In 2024 IEEE 63rd Conference on Decision and Control (CDC), pp. 5963-5968. IEEE, 2024.

[16] Shaked Regev, Nai-Yuan Chiang, Eric Darve, Cosmin G. Petra, Michael A. Saunders, Kasia Świrydowicz, and Slaven Peleš. HyKKT: a hybrid direct-iterative method for solving KKT linear systems. Optimization Methods & Software, pages 1–24, 2022.

[17] Świrydowicz, Kasia, Nicholson Koukpaizan, Maksudul Alam, Shaked Regev, Michael Saunders, and Slaven Peleš. “Iterative methods in GPU-resident linear solvers for nonlinear constrained optimization.” Parallel Computing 123 (2025): 103123.

[18] Kröger, Ole, Carleton Coffrin, Hassan Hijazi, and Harsha Nagarajan. "Juniper: An open-source nonlinear branch-and-bound solver in Julia." In International conference on the integration of constraint programming, artificial intelligence, and operations research, pp. 377-386. Cham: Springer International Publishing, 2018.

[19] Dunning, Iain, Joey Huchette, and Miles Lubin. "JuMP: A modeling language for mathematical optimization." SIAM review 59, no. 2 (2017): 295-320.